\documentclass[leqno,a4paper,titlepage,twoside,11pt,bibtex]{article}

\usepackage[T1]{fontenc}
\usepackage[latin1]{inputenc}
\usepackage[francais]{babel}

\usepackage{bbm,pifont,eufrak,eucal,mathrsfs,latexsym,graphicx,epsfig,color,rotating,multirow}

\usepackage{hyperref}

\usepackage[all]{xy}
\xyoption{2cell}

\newcounter{infra}[page]

\newenvironment{dem}[1][]{%
{\bf D\'emonstration #1 : }}{%
\hspace*{\fill}\nolinebreak[1]\hspace*{\fill}\underline{\bf Q.e.d.}\\}

\newenvironment{dem*}[1][]{%
{\bf D\'emonstration #1 : }}{%
 }

\newenvironment{eq*}{\begin{eqnarray*}}{\end{eqnarray*}}

\newenvironment{exemple}{{\sc Exemple : }}{}

\newtheorem{defi}{D\'efinition}

\newtheorem{thm}{Th\'eor\`eme}[section] \newtheorem{lem}[thm]{Lemme}
\newtheorem{pro}[thm]{Proposition}

\newcommand{\croi}{\times}

\newcommand{\A}{\ens A} 
\newcommand{\adh}{\overline}

\newcommand{\Aut}{\mathrm{Aut}}

 \newcommand{\boufl}{\xymatrix{ \!\!
\ar@(ur,dr) }}

\newcommand{\CCC}{{\adh{{\cal C}}}}
 
\newcommand{\cad}{\mbox{\it c-\`a-d }}
\newcommand{\call}{\mathscr}

\newcommand{\cf}{{\it cf. }}

 \newcommand{\cha}{\widehat}
 
\newcommand{\chap}{\widehat}

\newcommand{\Com}{\mathrm{Com}}

\newcommand{\cten}{\Box
\hspace*{-1.77ex}\raisebox{.19ex}{$\croi$}}

\newcommand{\dd}[1]{\frac{\partial}{\partial #1}}

\newcommand{\donne}{\mapsto}

\newcommand{\End}{\mathrm{End}} \newcommand{\ens}{\mathbbm}

 \newcommand{\et}{\mathrm{\; et \;}}

\newcommand{\GG}{\ens G}

\newcommand{\GGG}{\adh{\mathrm{PGL}_3}}

  \newcommand{\GL}{\mathrm{GL}} \newcommand{\goth}{\mathfrak}

\newcommand{\Id}{\mathrm{Id}}
 \newcommand{\ie}{\emph{i.e. }}

 \newcommand{\infi}{\infty}

\renewcommand{\int}{\mathrm{int\:}}

 \newcommand{\inv}{^{-1}}
 \newcommand{\iso}{\simeq}

\newcommand{\kk}{{\mathbf{k}}}

\newcommand{\ma}{\displaystyle} 

\newcommand{\mod}{\mathrm{ \: mod \:}}

\newcommand{\moins}{\:\setminus\:} 

 \newcommand{\Ou}{\mbox{\Ou}}

\newcommand{\p}{\:\:.}  
\newcommand{\PGL}{\mathrm{PGL}}
\newcommand{\PP}{\mathbbm{P}}

\newcommand{\Pic}{\mathrm{Pic}\;} 

\newcommand{\plus}{\oplus} \newcommand{\Plus}{\bigoplus}

\newcommand{\qq}{\forall\:}

\newcommand{\res}[1]{{\left | {}_{#1} \right.}}

\newcommand{\rg}{\mathrm{rg}\,}

\newcommand{\SL}{\mathrm{SL}}
\newcommand{\sll}{{\goth{sl}}}

\newcommand{\spec}{\mathrm{Spec}}  \newcommand{\sta}{\stackrel}
\newcommand{\sub}{\subseteq}

 \newcommand{\tens}{\otimes} 

\newcommand{\tenso}[1]{\raisebox{-1.5ex}{$\ma \stackrel{\displaystyle
\tens}{\scriptstyle #1}$}} 

\newcommand{\tilda}{\widetilde} 

\newcommand{\tq}{\: : \:}

\newcommand{\transposee}[1]{{\vphantom{#1}}^{\mathit t}{#1}}

\newcommand{\uni}{\cup}

\newcommand{\Vect}{\bigwedge\nolimits} 

\newcommand{\vect}{\land}

\newcommand{\ZZ}{\mbox{\Large $\mathbbm Z$}}
\newcommand{\zzz}{{\mathbf{z}}}

\begin{document}

\title{Repr\'esentations irr\'eductibles de certaines alg\`ebres d'op\'erateurs diff\'erentiels\\
Irreducible representations of some differential operator algebras}
 \author{Alexis TCHOUDJEM\\
  Institut Camille Jordan\\
Universit\'e Claude Bernard Lyon I\\
Boulevard du Onze Novembre 1918\\
69622 Villeurbanne\\
FRANCE\\
tchoudjem@math.univ-lyon1.fr
} \date{Villeurbanne, le \today}

\maketitle

{\bf Abstract :} For a projective variety $X$ and a line bundle $L$ over $X$, one considers the $L-$twisted global differential operator algebra $\call{D}_L(X)$ which naturally operates on the space of global sections $H^0(X,L)$. In the case where $X$ is the wonderful compactification of the group $\mathrm{PGL}_3$, one proves that the space $H^0(X,L)$ is an irreducible representation of the algebra $\call{D}_L(X)$ or zero. For that, one introduces a $2-$order differential operator which is defined over whole $X$ but which does not arise from the infinitesimal action of the automorphism group $\mathrm{Aut}(X)$.

\vskip 1cm

\section*{Introduction}

According to the famous Borel-Weil theorem, if $X$ is a flag variety and if $L$ is a line bundle over $X$, then the space $H^0(X,L)$ is an irreducible representation of the Lie algebra of the automorphism group of $X$, or $0$. If $\adh{G}$ is the wonderful compactification of an adjoint semisimple algebraic group $G$, then all the line bundles over $\adh{G}$ are $\cha{G} \croi \cha{G}-$linearized ($\cha{G}$ being the universal covering of $G$) and the space of global sections $H^0(\adh{G},L) $ is a representation of $\cha{G} \croi \cha{G}$ which is reducible in general. Now, the algebra $\call{D}_L(\adh{G})$ operates on $H^0(\adh{G},L)$ too ; that algebra is the global section algebra of a sheaf of algebras defined from $L$ and the sheaf of differential operators over $\adh{G}$. In the case where $G$ is $\PGL_3$, one gives a $2-$order differential operator defined over whole $\adh{G}$ ( \cf th\'eor\`eme \ref{thm:opdiffglob}). Then, by using this operator, one proves that the space $H^0(\adh{G},L)$ is irreducible (or zero) as a $\call{D}_L({\adh{G}})-$module, for all line bundles $L$ over $\adh{G}$. One finishes with a similar result over the \og complete conic variety \fg\ (\cf th\'eor\`eme \ref{thm:coniq}).

\vskip 1cm

{\bf R\'esum\'e : } \`A une vari\'et\'e projective $X$ et à un fibr\'e en droites $L$ sur $X$, on peut associer une alg\`ebre $\call{D}_L(X)$, l'alg\`ebre des $L-$op\'erateurs diff\'erentiels globaux sur $X$, qui op\`ere naturellement sur l'espace des sections globales $H^0(X,L)$. On montre ici que dans le cas particulier o\`u $X$ est la compactification magnifique du groupe  $\PGL_3$, l'espace $H^0(X,L)$ est soit nul soit une repr\'esentation irr\'eductible de l'alg\`ebre  $\call{D}_L(X)$. On introduit pour cela un op\'erateur diff\'erentiel d'ordre $2$ d\'efini sur tout $X$ mais qui ne provient pas de l'action infinit\'esimale sur $X$ du groupe d'automorphismes $\Aut X$.

\newpage

\tableofcontents

\vskip .5cm

Soit $ \kk$ un corps alg\'ebriquement clos de caract\'eristique nulle.

\addcontentsline{toc}{section}{Introduction}

\section*{Introduction}

D'apr\`es le c\'el\`ebre th\'eor\`eme de Borel-Weil, si $X$ est une vari\'et\'e de drapeaux et si $L$ est un fibr\'e en droites sur $X$, alors, lorsqu'il est non nul, l'espace $H^0(X,L)$ est une repr\'esentation irr\'eductible de l'alg\`ebre de Lie du groupe des automorphismes de $X$. Si $\adh{G}$ est la compactification magnifique d'un groupe alg\'ebrique semi-simple adjoint $G$, alors tout fibr\'e en droites $L$ sur $\adh{G}$ est $\cha{G} \croi \cha{G}-$lin\'earis\'e ($\cha{G}$ est le rev\^etement universel de $G$) et l'espace des sections globales $H^0(\adh{G},L)$ est une repr\'esentation de $\cha{G} \croi \cha{G}$ qui est r\'eductible en g\'en\'eral. Mais sur $H^0(\adh{G},L)$, l'anneau $\call{D}_L(\adh{G})$ op\`ere aussi ; c'est l'alg\`ebre des sections globales d'un faisceau d'alg\`ebres d\'efini à partir de $L$ et du faisceau des op\'erateurs diff\'erentiels sur $\adh{G}$. On donne ici dans le cas particulier o\`u $G = \PGL_3$ un op\'erateur diff\'erentiel d'ordre $2$ d\'efini sur tout $\adh{G}$ (\cf le th\'eor\`eme \ref{thm:opdiffglob}). Puis on d\'emontre en utilisant cet op\'erateur que l'espace $H^0(\adh{G},L)$ est irr\'eductible (ou nul) comme $\call{D}_L(\adh{G})-$module pour tout fibr\'e en droites $L$ sur $\adh{G}$. On termine par un r\'esultat analogue sur la vari\'et\'e des \og coniques compl\`etes \fg\ (\cf le th\'eor\`eme \ref{thm:coniq}).

\section{Op\'erateurs diff\'erentiels sur une vari\'et\'e}

Soit $X$ une vari\'et\'e alg\'ebrique lisse.

On rappelle ici la d\'efinition du faisceau $\call{D}_X$ des op\'erateurs diff\'erentiels sur $X$. 

On note $\call{E}nd_\kk(\call{O}_X)$ le faisceau des $\kk-$endomorphismes locaux  de $\call{O}_X$. On identifiera $\call{O}_X$ \`a un sous-faisceau $\call{E}nd_\kk(\call{O}_X)$.

On pose : $\call{D}_X^{(0)}:=\call{O}_X$ et pour tout $n \ge 0$ et tout ouvert $U$ de $X$ :
\[\Gamma(U,\call{D}_X^{(n+1)}):= \]

\[ \left\{ d \in \End_\kk(\call{O}_U) \tq \qq V \sub U \mbox{ ouvert }, \qq a \in \call{O}_X(U),\, [d\res{V},a] \in \call{D}_X^{(n)}(V)\right\} \p\]

\begin{defi}
\[\call{D}_X:=\uni_{n \ge 0} \call{D}_X^{(n)} \p\]
\end{defi}

On notera $\call{D}(X)$ l'anneau $\Gamma(X,\call{D}_X)$ des op\'erateurs diff\'erentiels globaux sur $X$.

\begin{exemple}
Si $\PP^1 = \{[x_0:x_1] \tq (x_0,x_1) \in \kk^2 \moins\{(0,0)\} \}$, alors :
\[{\call{D}(\PP^1)} = \kk[x_0\partial_{x_0} , x_0\partial_{x_1}, x_1 \partial{x_0}]\]

(avec $x_1\partial_{x_1} = -x_0\partial_{x_0}$).
\end{exemple}

\subsection{Op\'erateurs diff\'erentiels tordus par un faisceau inversible}

Soit $\call{L}$ un faisceau inversible sur $X$. On rappelle que le faisceau d'algèbres des op\'erateurs diff\'erentiels tordus par $\call{L}$, not\'e $\call{D}_{\call{L}}$, est un sous-faisceau de $\call{E}nd_\kk(\call{L})$ d\'efini comme le faisceau $\call{D}_X$, par r\'ecurrence :

On pose $\call{D}_{\call{L}}^{(0)} := \call{O}_X$ et pour tout $n \ge 0$ et tout $U$ ouvert de $X$ :
\[\Gamma(U,\call{D}_\call{L}^{(n+1)}):= \] \[ \left\{d \in \End_\kk(\call{L}\res{U}) \tq \qq V \sub U \mbox{ ouvert },\, \qq a \in \call{O}_X(V),\, [d\res{V},a]\in \call{D}_\call{L}^{(n)}(V) \right\} \p\]

\begin{defi}
\[\call{D}_\call{L} := \uni_{n\ge 0}\call{D}_\call{L}^{(n)} \p\]
\end{defi}

\begin{exemple}
Soient $n \in \ZZ$ et $\call{L}_n:=\call{O}_{\PP^1}(n)$ le faisceau inversible des fonctions $n-$homog\`enes. On a :
\[\call{D}_{\call{L}_n}(\PP^1) = \kk[x_0\partial_{x_0} , x_0\partial_{x_1}, x_1 \partial{x_0}]\]
(avec $x_1\partial_{x_1} = n -x_0\partial_{x_0}$).

\end{exemple}

{\bf Remarques :} 

1) On a un isomorphisme de faisceaux d'algèbres : \[\call{D}_\call{L} \iso \call{L}\tenso{\call{O}_X} \call{D}_X \tenso{\call{O}_X} \call{L}\inv \p\]

2) L'algèbre $\call{D}_\call{L}(X)$ opère naturellement sur tous les groupes de cohomologie $H^i(U,\call{L})$, pour tout $i \ge 0$ et tout $U$ ouvert de $X$.

\subsection{Action de groupes}

 Si $G$ est un groupe alg\'ebrique qui agit alg\'ebriquement sur $X$ et si $\call{L}$ est un faisceau inversible $G-$lin\'earis\'e sur $X$, alors l'action infinit\'esimale de l'alg\`ebre de Lie de $G$, $\goth g$,  sur $X$ (\cf \cite{kempf}[\S 1]) induit un morphisme naturel d'algèbres associatives :
\[U(\goth g) \to \call{D}_\call{L}(X)\]
(o\`u $U(\goth g)$ est l'algèbre enveloppante de $\goth g$). 

\begin{exemple}
--- Si $\ma X = \PP^1$ et $G = SL_2$, $\sll_2$ son alg\`ebre de Lie, alors on a un morphisme surjectif :

\[U(\sll_2) \to \call{D}_{\call{L}_n}(\PP^1)\]
\[ \scriptsize{\left(\begin{array}{cc}
0&1\\
0&0
\end{array}\right) \, , \, \left(\begin{array}{cc}
0&0\\
1&0
\end{array}\right) \, , \, \left(\begin{array}{cc}
1&0\\
0&-1
\end{array}\right) \donne -x_1\partial_{x_0} , -x_0 \partial_{x_1}, -x_0 \partial_{x_0} +x_1\partial_{x_1} \p}\]

--- En revanche, si $X$ est l'\'eclatement de $\PP^2$ en $0:=[1:0:0]$, si $G:=\Aut X = \scriptsize{\bigg \{ \left[\begin{array}{ccc}
* & * &*\\
0 &* &*\\
0 & * & *
\end{array}\right] \in \PGL_3(\kk) \bigg\}}$, le morphisme :

\[U(\goth g) \to \call{D}(X) \]
n'est pas surjectif.
\end{exemple}

\section{Compactifications magnifiques des groupes adjoints}

Soit $G$ un groupe alg\'ebrique lin\'eaire semi-simple et adjoint. On notera $R :\chap{G} \to G$ son rev\^etement universel ($\chap{G}$ est semi-simple simplement connexe et $\ker R$ est le centre de $\chap{G}$). On notera $\goth g$ l'alg\`ebre de Lie de $G$.

Les multiplications \`a gauche et \`a droite d\'efinissent une action de $G \croi G$ sur $G$. D'apr\`es \cite[\S 2]{dcp}, il existe une compactification magnifique de $G$ \cad une vari\'et\'e projective, lisse et connexe, not\'ee $\adh{G}$ qui contient $G$ telle que :

i) $G \croi G$ agit sur $\adh{G}$ et cette action prolonge l'action sur $G$ ;

ii) $G$ est ouvert dans $\adh{G}$ ;

iii) si on note $Z_1,...,Z_r$ les composantes irr\'eductibles de $\adh{G} \moins G$ (les {\it diviseurs limitrophes}), alors $\adh{G}\moins G = Z_1 \uni .. \uni Z_r$ est un diviseur \`a croisements normaux et $r $ est le rang de $G$ ;

iv) chaque adh\'erence de $G \croi G-$orbite de $\adh{G}$ est l'intersection (transverse) des $Z_i$ qui la contiennent ;

v) l'intersection $ Z_1 \cap ... \cap Z_r$ est l'unique orbite ferm\'ee de $\adh{G}$.

\section{Immersion dans un produit d'espasces projectifs}

Soient $\chap{T} \sub \chap{B}$ un tore maximal et un sous-groupe de Borel de $\chap{G}$. On pose $T:=R(\chap{T}), B:=R(\chap{B})$. Soient $\alpha_1, ..., \alpha_r$ la base correspondante du syst\`eme de racines de $(G,T)$ (les $\alpha_i$ sont en fait des caract\`eres de $T$). Soient $\omega_1,...,\omega_r$ les poids fondamentaux correspondant \`a cette base (ce sont des caract\`eres de $\chap{T}$). Notons $\rho_1 : \chap{G} \to \GL(V_{\omega_1}), ..., \rho_r : \chap{G} \to \GL( V_{\omega_r})$ les repr\'esentations irr\'eductibles de $\chap{G}$ associ\'ees.

Posons aussi pour tout $i : E_{\omega_i} := \End(V_{\omega_i})$.

Consid\'erons \[i : G \to  \PP(E_{\omega_1}) \croi ... \croi \PP(E_{\omega_r})\]
\[g \donne ([\rho_1(g)] ,..., [\rho_r(g)]) \p\]

Le morphisme $i$ est une immersion et on peut identifier  $\adh{G}$ à $\adh{i(G)}$.

On notera pour tous $g,g'\in G$, $x \in \adh{G}$, $gxg':=(g,{g'}\inv).x$.

\subsection{Notations}
Fixons quelques notations.

Soient $\chap{B}^-$ (resp. $B^-$) le sous-groupe de Borel oppos\'e \`a $\chap{B}$ relativement \`a $\chap{T}$ (resp. oppos\'e \`a $B$ relativement \`a $T$).

On notera aussi $U$ et $U^-$ les radicaux unipotents de $B$ et $B^-$,  $\goth n$ et $\goth n^-$ les alg\`ebres de Lie de $U$ et $U^-$.

La vari\'et\'e $\adh{G}$ poss\`ede un unique point fixe pour $B^- \croi B$. Notons le $\zzz$. Pour tout caract\`ere $\lambda$ de $\chap{T}$ dominant (par rapport \`a $\chap{B}$), on note $V_\lambda$ le $\chap{G}-$module irr\'eductible de plus haut poids $\lambda$ et $\lambda^*$ le plus haut poids du $\chap{G}-$module irr\'eductible dual $V_\lambda^*$. On choisit $v_\lambda \in V_\lambda$ et $v^*_{-\lambda} \in V_\lambda^*$ un $\chap{B}-$vecteur propre de poids $\lambda$ et un $\chap{B}^--$vecteur propre de poids $-\lambda$ tels que $\langle v^*_{-\lambda},v_\lambda \rangle = 1$. De m\^eme, on choisit un $\chap{B}^-$vecteur propre $v_{-\lambda^*} \in V_\lambda$ et un $\chap{B}-$vecteur propre $v_{\lambda^*}^*\in V_\lambda^*$. Si on consid\`ere $\adh{G}$ comme une sous-vari\'et\'e ferm\'ee de $\PP(E_{\omega_1}) \croi .. \croi \PP(E_{\omega_r})$ et compte tenu de l'identification : 
\[E_{\omega_i} = V_{\omega_i}
\tens V_{\omega_i}^*\]
on a : 
\[\zzz = ([v_{-\omega_1^*} \tens v^*_{\omega_1^*}],...,[v_{-\omega_r^*} \tens v^*_{\omega_r^*}]) \p\]

\section{La grosse cellule}\label{sec:grocel}

 Il existe un ouvert (unique) de $\adh{G}$, not\'e $(\adh{G})_0$ et appel\'e la {\it grosse cellule}, qui est $B \croi B^-$ invariant et isomorphe \`a un espace affine.

De plus, $\adh{G} = (G\croi G).(\adh{G})_0$ et $BB^-  = (\adh{G})_0 \cap G$.

De plus, la d\'ecomposition $BB^-= UTU^- \iso U \croi T \croi U^-$ se \og prolonge \fg\ \`a $(\adh{G})_0$. En effet, si on pose $(\adh{T})_0 := \adh{T} \cap (\adh{G})_0$, on a un isomorphisme :

\[U \croi (\adh{T})_0 \croi U^- \to (\adh{G})_0\]

\[(u,x,u') \donne uxu' \p\]

Ajoutons que pour tout $1 \le i \le r$, le caract\`ere $\alpha_i : T \to \GG_m$ se prolonge en un morphisme, encore not\'e $\alpha_i : (\adh{T})_0 \to \A^1$ et que l'on a aussi un isomorphisme :
\[(\adh{T})_0 \to \A^r\]
\[x \donne (\alpha_1(x),...,\alpha_r(x) )\p\]

Pour simplifier, on notera encore $\alpha_i$ les morphismes :
\[(\adh{G})_0 \to \A^1\]
\[uxu' \donne \alpha_i(x)\]

(pour tout $u \in U$ et tout $u'\in U^-$). 

En particulier, $\kk [(\adh{G})_0] = \kk [U \croi U^-][\alpha_1,...,\alpha_r]$.

En fait, on a aussi :
\[(\adh{G})_0 = \left\{ ([A_1],...,[A_r]) \in \adh{G} \sub \PP(E_{\omega_1}) \croi ... \croi \PP(E_{\omega_r}) \tq \qq 1 \le i \le r, f_i(A_i) \neq 0 \right\}\p\]

o\`u pour tout $i$, $f_i$ est la forme lin\'eaire :
\[f_i : E_{\omega_i} \to \kk\;, \; A \donne \langle v^*_{\omega_i^*}
, A v_{-\omega_i^*}\rangle \p\]

De plus, les diviseurs limitrophes $Z_i$ v\'erifient :

\[Z_i \cap (\adh{G})_0 = \left\{ x \in (\adh{G})_0 \tq \alpha_i(x) =0 \right\} \p\]

\section{Le groupe de Picard}
Rappelons la description des fibr\'es en droites (ou des faisceaux inversibles) sur $\adh{G}$.

Tous les faisceaux inversibles sur $\adh{G}$ sont $\chap{G}\croi \chap{G}-$lin\'earis\'es de mani\`ere unique de plus :

\begin{pro}[{\cf \cite[\S 8]{dcp}}]
On a un isomorphisme :

\[\ZZ^r \to \Pic(\adh{G})\]
\[(n_1,...,n_r) \donne \left( \call{O}_{\PP(E_{\omega_1})}(n_1) \cten ... \cten \call{O}_{\PP(E_{\omega_r})}(n_r)\right) \res{\adh{G}} \p\]
\end{pro}

Si $\lambda = n_1 \omega_1 + ... + n_r\omega_r$ est un caract\`ere de $\chap{T}$, on posera :
\[\call{L}_\lambda := \left( \call{O}_{\PP(E_{\omega_1})} (n_1)\cten ... \cten \call{O}_{\PP(E_{\omega_r})}(n_r)\right) \res{\adh{G}} \p\]

En particulier, sur la fibre $\call{L}_\lambda\res{\zzz}$, le tore $\chap{T}\croi \chap{T}$ agit avec le poids $(\lambda^*,-\lambda^*)$ o\`u si l'on note $W$ le groupe de Weyl de $(G,T)$ et $w_0$ l'\'el\'ement le plus long de $W$, $\lambda^*:=-w_0\lambda$.

{\bf Remarque :} si on note $\call{O}_{\adh{G}}(Z_i)$ le faisceau inversible associ\'e au diviseur $Z_i$, on a un isomorphisme de faisceaux $G \croi G-$lin\'earis\'es : \[\call{O}_{\adh{G}}(Z_i) \iso \call{L}_{\alpha_i^*} \p\]

\section{Sections des faisceaux inversibles} 

\subsection{sur la grosse cellule}\label{par:secglob}
La grosse cellule $(\adh{G})_0$ est isomorphe \`a un espace affine donc pour tout caract\`ere $\lambda$ de $\chap{T}$, $\Gamma((\adh{G})_0, \call{L}_\lambda)$ est un $k[(\adh{G})_0]-$module libre de rang $1$. Si $\lambda= n_1 \omega_1 + ... + n_r \omega_r$, alors :
\[\Gamma((\adh{G})_0, \call{L}_\lambda) = \kk [(\adh{G})_0]f_\lambda\]
o\`u  $f_\lambda :=(f_1^{n_1} \cten ... \cten f_r^{n_r} ) \res{(\adh{G})_0}$.

\subsection{sections globales}

Rappelons que l'on consid\`ere les caract\`eres  $\alpha_i$,  $(1 \le i \le r)$, comme des fonctions r\'eguli\`eres $U \croi U^--$invariantes sur la grosse cellule $(\adh{G})_0$ (\cf \S \ref{sec:grocel}) .
\begin{thm}[{\cf \cite[\S 8.2]{dcp}}]
Pour tout caract\`ere $\lambda$ de $\chap{T}$, on a un isomorphisme de $U(\goth g \croi \goth g)-$modules :
\[
\Gamma(\adh{G}, \call{L}_\lambda) \iso \Plus_{m_1 ,...,m_r \ge 0\atop \lambda^* - m_1 \alpha_1 - ... - m_r \alpha_r \mathrm{\, dominant}} U(\goth g \croi \goth g) \alpha_1^{m_1}...\alpha_r^{m_r} f_\lambda \p\]
\end{thm}

\section{Cas particulier}

On suppose maintenant que $G = \PGL_3(\kk)$. Dans ce cas :

\[\chap{G}= \SL_3(\kk), r=2\]
Les sous-groupes $\chap{T},T,\chap{B},B,\chap{B}^-,B^-$ sont respectivement les sous-groupes des matrices diagonales, triangulaires sup\'erieures et inf\'erieures. Si on note $\epsilon_i$ le caract\`ere :

\[\left(\begin{array}{ccc}
t_1 & &\\
&t_2&\\
&&t_3
\end{array}\right) \donne t_i\]
on a : $\omega_1 = \epsilon_1,\omega_2=\epsilon_1+\epsilon_2, \alpha_1 = \epsilon_1 -\epsilon_2,\alpha_2 = \epsilon_2 -\epsilon_3$.

On a avec ces notations :
 \[E_{\omega_1} = \End(\kk^3) \et E_{\omega_2}=\End(\Vect^2 \kk^3) \iso \End(\kk^3 )\p\]
Donc si l'on prend pour base de $\kk^3$ la base canonique $e_1,e_2,e_3$ et pour base de $\Vect^2 \kk^3$ la base canonique duale : $e_2 \vect e_3,e_3 \vect e_1, e_1 \vect e_2$, on a $E_{\omega_1} \iso E_{\omega_2} \iso M_3(\kk)$.

En particulier, le plongement de $G$ dans $\adh{G}$ est donn\'e par :
\[G \to \PP^8(\kk) \croi \PP^8(\kk)\]
\[[g] \donne ([g],[\Com(g)])\] 
o\`u $\Com(g)$ est la comatrice de $g$ .

On en d\'eduit une description explicite de $\adh{\PGL_3(\kk)}$ avec des \'equations :

\begin{pro}
Si $G =\PGL_3(\kk)$, alors :
\[\adh{G} = \left\{ ([g],[g']) \in \PP(M_3(\kk)) \croi \PP(M_3(\kk)) \tq g{}^t\!g' \in \kk I_3 \right\} \]
o\`u $I_3$ est la matrice identit\'e $3 \croi 3$.
\end{pro}

De plus, si on note $g_{i,j}$ le $(i,j)-$i\`eme coefficient d'une matrice $g$, on a :

\[(\adh{G})_0 = \left\{([g],[g']) \in \adh{G} \tq g_{3,3}g'_{1,1} \neq 0\right\}\]
\[(\adh{T})_0 = \left\{{\left({\left[\begin{array}{ccc}
t_2t'_2&&\\
&t_2&\\
&&1

\end{array}\right] ,\left[\begin{array}{ccc}
1 &&\\
&t'_2&\\
&&t_2t'_2

\end{array}\right]}\right)}\tq t_2,t'_2 \in \kk \right\} \p\]

Soit $a:=\left({\left[\begin{array}{ccc}
t_1&&\\
&t_2&\\
&&t_3
\end{array}\right],\left[\begin{array}{ccc}
t'_1&&\\
&t'_2&\\
&&t'_3
\end{array}\right]}\right)
 \in (\adh{T})_0$. Pour tout \[
u=\left[\begin{array}{ccc}1&u_{1,2}&u_{1,3}\\
&1&u_{2,3}\\
&&1
\end{array}\right] \in U \;
,\mbox{ tout } u'=\left[\begin{array}{ccc}
1&&\\
u_{2,1}&1&\\
u_{3,1}&u_{3,2}&1
\end{array}\right]\in  U^-\; ,\]
on a : 
\[\alpha_1(uau') = \frac{t'_2}{t'_1}, \alpha_2(uau') =\frac{t_2}{t_3}\]
et on pose $U_{i,j}(uau') := u_{i,j}$ si $1 \le i\not=j \le 3$.

Ainsi, $\kk[(\adh{G})_0]= k[U_{i,j} \tq 1\le i\neq j \le 3][\alpha_1,\alpha_2]$.

On remarque enfin que dans ce cas  particulier, les diviseurs $Z_1$ et $Z_2$ sont d\'efinis par :

\[Z_1 = \left\{([g],[g']) \in \adh{G} \tq \rg (g) =1 \right\} \; ,\; Z_2 = \left\{([g],[g'] \in \adh{G} \tq \rg (g')=1 \right\} \p\]
\section{Formules de changement de variables}

\subsection{Fonctions}
Puisque \[G \cap (\adh{G})_0
 = BB^-=\left\{[g] \in G \tq g_{3,3}\neq 0, g_{2,2}g_{3,3}-g_{2,3}g_{3,2} \neq 0\right\} \]
les fonctions $\alpha_1,\alpha_2, U_{i,j}, 1 \le i\neq j \le 3$ sont r\'eguli\`eres sur $BB^-$ On peut donc les exprimer en fonction des coordonn\'ees $g_{i,j}$ de la matrice $g$. En posant $\Delta_{i,j}(g):=(-1)^{i+j}\Com(g)_{i,j}$ et $\Delta(g):=\det (g)$, on trouve :
\begin{pro}\label{pro:cdv}
Pour tout $g \in BB^{-}$, \[\alpha_1(g) = \frac{g_{3,3}\Delta(g)}{\Delta_{1,1}^2(g)}\]
\[\alpha_2(g) = \frac{\Delta_{1,1}(g)}{g_{3,3}^2}\]
\[U_{1,2}(g)= \frac{\Delta_{2,1}(g)}{\Delta_{1,1}(g)}\]

\[U_{2,1}(g)= \frac{\Delta_{1,2}(g)}{\Delta_{1,1}(g)}\]
\[U_{1,3}(g)=\frac{g_{1,3}}{g_{3,3}}\]
\[U_{3,1}(g)=\frac{g_{3,1}}{g_{3,3}}\]
\[U_{2,3}(g)=\frac{g_{2,3}}{g_{3,3}}\]
\[U_{3,2}(g)=\frac{g_{3,2}}{g_{3,3}} \]

et r\'eciproquement :

\[\frac{g_{1,1}}{g_{3,3}}= \alpha_1\alpha_2 + \alpha_2 U_{1,2}U_{2,1}+U_{1,3}U_{3,1}\]
\[\frac{g_{1,2}}{g_{3,3}}=\alpha_2U_{1,2}+U_{1,3}U_{3,2}\]
\[\frac{g_{1,3}}{g_{3,3}}=U_{1,3}\]
\[\frac{g_{2,1}}{g_{3,3}}=\alpha_2 U_{2,1}+U_{2,3}U_{3,1}\]
\[\frac{g_{2,2}}{g_{3,3}}=\alpha_2+U_{2,3}U_{3,2}\]
\[\frac{g_{2,3}}{g_{3,3}}=U_{2,3}\]
\[\frac{g_{3,1}}{g_{3,3}}=U_{3,1}\]
\[\frac{g_{3,2}}{g_{3,3}}=U_{2,3}\]

\end{pro}

\subsection{D\'erivations issues de l'action de l'alg\`ebre de Lie}

On choisit la base suivante de $\goth g = \sll_3$ :

\[X_1:= \left(\begin{array}{ccc}
0&1&0\\
0&0&0\\
0&0&0
\end{array}\right) \; ,\;  X_2:= \left(\begin{array}{ccc}
0&0&0\\
0&0&1\\
0&0&0
\end{array}\right) \;,\; X_3:= \left(\begin{array}{ccc}
0&0&1\\
0&0&0\\
0&0&0
\end{array}\right)\]
\[Y_1:=  \left(\begin{array}{ccc}
0&0&0\\
1&0&0\\
0&0&0
\end{array}\right)\; ,\;  Y_2 := \left(\begin{array}{ccc}
0&0&0\\
0&0&0\\
0&1&0
\end{array}\right) \; , \; Y_3:=  \left(\begin{array}{ccc}
0&0&0\\
0&0&0\\
1&0&0
\end{array}\right)\]
\[H_1:=  \left(\begin{array}{ccc}
1&0&0\\
0&-1&0\\
0&0&0
\end{array}\right) \; ,\; H_2:=  \left(\begin{array}{ccc}
0&0&0\\
0&1&0\\
0&0&-1
\end{array}\right)
\]

Pour tout $\xi \in U(\goth g)$, on note $\xi^{(g)}$ l'image de $(\xi,0) \in U(\goth g \croi 0)$ dans $U(\goth g \croi \goth g)$. Soit $\Phi_0 $ le morphisme d'alg\`ebres :
\[\Phi_0 : U(\goth g \croi \goth g) \to \call{D}(\adh{G})\]
induit par l'action de $\goth g \croi \goth g$ sur $\call{O}_{\adh{G}}$.

Comme $\call{D}({\adh{G}}) \sub \Gamma ({G},\call{D}_{\adh{G}}) \cap \Gamma((\adh{G})_0,\call{D}_{\adh{G}})$, on peut exprimer les $\Phi_0(\xi^{(g)})$, $\xi \in \goth g$, en fonction des coordonn\'ees de $G$ et en fonction des coordonn\'ees de $(\adh{G})_0$ :

\begin{pro}
Dans l'anneau $ \Gamma ({G},\call{D}_{\adh{G}})$, on a :

\begin{eqnarray*}
\Phi_0(Y_1^{(g)}) & = & -g_{1,1}\dd{g_{2,1}} -g_{1,2}\dd{g_{2,2}} -g_{1,3}\dd{g_{2,3}}\;,\\
\Phi_0(Y_2^{(g)}) & = & -g_{2,1}\dd{g_{3,1}} -g_{2,2}\dd{g_{3,2}} -g_{2,3}\dd{g_{3,3}} \;,\\
\Phi_0(Y_3^{(g)}) & = & -g_{3,1}\dd{g_{2,1}} -g_{1,2}\dd{g_{3,2}} -g_{1,3}\dd{g_{3,3}}\;,\\
\Phi_0(X_1^{(g)}) & = & -g_{2,1}\dd{g_{1,1}} -g_{2,2}\dd{g_{1,2}} -g_{2,3}\dd{g_{1,3}}\;,\\
\Phi_0(X_2^{(g)}) & = & -g_{3,1}\dd{g_{2,1}} -g_{3,2}\dd{g_{2,2}} -g_{3,3}\dd{g_{2,3}}\;,\\
\Phi_0(X_3^{(g)}) & = & -g_{3,1}\dd{g_{1,1}} -g_{3,2}\dd{g_{1,2}} -g_{3,3}\dd{g_{1,3}} \;,\\
\Phi_0(H_1^{(g)}) & = & -g_{1,1}\dd{g_{1,1}} - g_{1,2}\dd{g_{1,2}} - g_{1,3}\dd{g_{1,3}} + g_{2,1}\dd{g_{2,1}} +g_{2,2}\dd{g_{2,2}} + g_{2,3}\dd{g_{2,3}}\;,\\
\Phi_0(H_2^{(g)}) & = & - g_{2,1}\dd{g_{2,1}} - g_{2,2}\dd{g_{2,2}} - g_{2,3}\dd{g_{2,3}} + g_{3,1}\dd{g_{3,1}} +g_{3,2}\dd{g_{3,2}} + g_{3,3}\dd{g_{3,3}} \p
\end{eqnarray*}

Dans l'anneau $ \Gamma ((\adh{G})_0,\call{D}_{\adh{G}})$, on a :

\begin{eqnarray*} 
{ \Phi_0(Y_1^{(g)})} & =& 2U_{1,2}\alpha_1\dd{\alpha_1} -U_{1,2}\alpha_2\dd{\alpha_2}\\
&& -U_{1,3}\dd{U_{2,3}} -U_{1,2}^2\dd{U_{1,2}} \\
&&-\alpha_1\dd{U_{2,1}}  \;,\\
{ \Phi_0(Y_2^{(g)})}& =& -U_{2,3} \alpha_1 \dd{\alpha_1} + 2U_{2,3} \alpha_2 \dd{\alpha_2} \\
&& + (U_{1,3}-U_{1,2}U_{2,3})\dd{U_{1,2}} +U_{2,3}^2\dd{U_{2,3}} - U_{1,3}U_{2,3} \dd{U_{1,3}} \\
&&- \alpha_2 ( U_{2,1}\dd{U_{3,1}} +\dd{U_{3,2}}) \;,\\
{\Phi(Y_3^{(g)})}
& =& (U_{1,3}-2U_{1,3}U_{3,2}) \alpha_1\dd{\alpha_1} + (U_{1,3}+U_{1,2}U_{2,3})\alpha_2\dd{\alpha_2} \\&&+ (U_{1,3}-U_{1,2}U_{2,3})U_{1,2}\dd{U_{1,2}} + U_{1,3}U_{2,3}\dd{U_{2,3}} + U_{1,3}^2\dd{U_{1,3}} \\
&&+\alpha_1 U_{2,3}\dd{U_{2,1}} - \alpha_2(U_{1,2} \dd{U_{3,2}} + U_{1,2}U_{2,1}\dd{U_{3,1}}) - \alpha_1 \alpha_2 \dd{U_{3,1}} \;,\\
\Phi_0(X_1^{(g)}) &= &-\dd{U_{1,2}} - U_{2,3}\dd{U_{1,3}} \;,\\
\Phi_0(X_2^{(g)}) &=& -\dd{U_{2,3}} \;,\\
\Phi_0(X_3^{(g)}) &=& -\dd{U_{1,3}} \p 
\end{eqnarray*}

\end{pro}

\begin{dem}
Par exemple, pour la deuxi\`eme liste de formules, on utilise que :
\[\Phi_0(\xi^{(g)}) = (\xi^{(g)} . \alpha_1 )\dd{\alpha_1} + (\xi^{(g)} . \alpha_2) \dd{\alpha_2} + \sum_{1 \le i \neq j \le n}(\xi^{(g)}.U_{i,j} )\dd{U_{i,j}} \]
pour tout $\xi \in \goth g$.
\end{dem}

\subsection{Un op\'erateur diff\'erentiel d'ordre $2$}

Comme $(\adh{G})_0$ est isomorphe \`a un espace affine,
\[\Gamma ( (\adh{G})_0,\call{D}_{\adh{G}}) = \kk[\alpha_1,\alpha_2, U_{i,j} \tq 1 \le i\neq j \le 3][\partial_{\alpha_1},\partial_{\alpha_2},\partial_{U_{i,j}} \tq  1 \le i\neq j \le 3] \p\]

De plus, comme $BB^-$ est un ouvert de $G$,
\[\Gamma((\adh{G})_0,\call{D}_{\adh{G}}) \sub \Gamma(BB^-,\call{D}_{\adh{G}}) \sub \kk(g_{i,j} \tq 1 \le i ,j \le 3)[\partial_{g_{k,l}}\tq 1 \le k,l \le 3] \p
\]

On peut donc exprimer en particulier $\partial_{\alpha_1}$ et $\partial_{\alpha_2}$ en fonction des fonctions coordonn\'ees $g_{i,j}$ et des d\'eriv\'ees partielles $\partial_{g_{k,l}}$ :

\begin{pro}
On a :
\[\partial_{\alpha_1} = \frac{\Delta_{1,1}}{g_{3,3}}\partial_{g_{1,1}}\]
\[\partial_{\alpha_2} = \frac{g_{3,3}}{\Delta_{1,1}}\left( \Delta_{2,2}\partial_{g_{1,1}} + \Delta_{1,1}\partial_{g_{2,2}} + \Delta_{2,1}\partial_{g_{1,2}} + \Delta_{1,2}\partial_{g_{2,1}} \right) \p\]

\end{pro}

\begin{dem}
Posons $(F_1,...,F_8) := (\alpha_1,\alpha_2, U_{1,2},...,U_{3,2})$ 

et $(x_1,...,x_8) := (\frac{g_{1,1}}{g_{3,3}}, ... , \frac{g_{3,2}}{g_{3,3}})$. Les coefficients de l'inverse de la matrice jacobienne :
\[\left(\frac{\partial F_i }{\partial x_j} \right)_{1 \le i,j \le 8}\]

donnent les coefficients des d\'erivations $\partial\alpha_1$, ... en fonction des d\'erivations $\partial_{\frac{g_{i,j}}{g_{3,3}}} = g_{3,3}\partial_{g_{i,j}}$, $1 \le i ,j \le 3$, $(i,j) \neq (3,3)$.
\end{dem}
On remarque en particulier que :
\[\partial_{\alpha_1}\partial_{\alpha_2} = \partial_{g_{1,1}}\left( \Delta_{2,2}\partial_{g_{1,1}} + \Delta_{1,1}\partial_{g_{2,2}} + \Delta_{2,1}\partial_{g_{1,2}} + \Delta_{1,2}\partial_{g_{2,1}} \right) \in \Gamma(G, \call{D}_{\adh{G}}) \p\]

On en d\'eduit :

\begin{thm}\label{thm:opdiffglob}
L'op\'erateur diff\'erentiel \[D_0:=\partial_{\alpha_1}\partial_{\alpha_2} = \partial_{g_{1,1}}\left( \Delta_{2,2}\partial_{g_{1,1}} + \Delta_{1,1}\partial_{g_{2,2}} + \Delta_{2,1}\partial_{g_{1,2}} + \Delta_{1,2}\partial_{g_{2,1}} \right)\]
est d\'efini sur $\adh{G}$ tout entier.
\end{thm}

{\bf Remarques :}

--- Les op\'erateurs diff\'erentiels $\partial_{\alpha_i}$ ($i=1$ ou $2$) ne sont pas dans $\call{D}({\adh{G}})$ ;

--- on peut v\'erifier que l'op\'erateur diff\'erentiel $D_0$ n'est pas l'image d'un \'el\'ement de $U(\goth g \croi \goth g)$ car, par exemple, il ne commute pas avec l'\'el\'ement de Casimir standard de $U(\goth g \croi 0)$ (ou $U(0 \croi \goth g)$). 

\begin{dem}
Puisque $D_0$ appartient au $G\croi G-$module rationnel $ \Gamma(G,\call{D}_{\adh{G}})$, $D_0$ est $U(\goth n^- \croi \goth n^-)-$fini et $U(\goth n\croi \goth n)-$fini. Donc d'apr\`es le lemme de prolongement \ref{lem:prolonge} qui suit, il existe un ouvert $\Omega$, contenant $(\adh{G})_0$, $U^- \croi U^-$ et $U\croi U-$stable, tel que $D_0 \in \Gamma(\Omega,  \call{D}_{\adh{G}})$. Or, comme $G$ est engendr\'e par $U$ et $U^-$, $\Omega$ est $G \croi G-$stable et $\Omega = (G \croi G).(\adh{G})_0 = \adh{G}$. 

\end{dem}

Pour le lemme suivant, on note $\goth g_a$ l'alg\`ebre de Lie du groupe $\GG_a$ et pour tout $i \ge 0$, $\goth g_a^i$le sous-espace $\kk d^i$ de l'alg\`ebre enveloppante $U(\goth g_a)$ (pour un g\'en\'erateur quelconque $d$ de $\goth g_a$).

\begin{lem}[de prolongement]\label{lem:prolonge}
Soit $Y$ une $\GG_a-$vari\'et\'e. Soit $\call{F}$ un faisceau localement libre de $\call{O}_Y-$modules et $\GG_a-$lin\'earis\'e sur $Y$. Soit $\Omega$ un ouvert de $Y$. Alors la $\GG_a-$lin\'earisation de $\call{F}$ induit une structure de $U(\goth g_a)-$module sur  $\Gamma(\Omega,\call{F})$. 

Soit $\sigma \in \Gamma(\Omega, \call{F})$. Si pour un certain $n_0 \ge 0$  $\goth g_a^{n_0}.\sigma = 0$, alors il existe $\tilda{\Omega}$ un ouvert $\GG_a-$stable de $Y$ et $\tilda{\sigma} \in \Gamma(\tilda{\Omega},\call{F})$ tel que :

\[\tilda{\sigma}\res{\Omega } = \sigma\]
(autrement dit $\sigma$ se prolonge en une section d\'efinie sur un ouvert $\GG_a-$stable). 
\end{lem}

\begin{dem}
Notons $\mu : \GG_a \croi Y \to Y$ le morphisme d\'efini par l'action : $(g,y) \donne g.y$ et $p : \GG_a \croi Y \to Y$ la projection sur $Y$. Le faisceau $\call{F}$ est $\GG_a-$lin\'earis\'e : cela signifie qu'il existe un isomorphisme de faisceaux de $\call{O}_{\GG_a\croi Y}-$modules :

\[\Phi : \mu^* \call{F} \sta{\iso}{\to} p^*\call{F} \]
qui v\'erifie les conditions de \cite[1 \S 3 def. 1.6]{Mum}.

Pour tout $n > 0$, soit $\GG_{a,n} := \spec\left(\kk[T]/(T^n)\right)$ le $n-$i\`eme voisinage infinit\'esimal de $0$ dans $\GG_a$.La restriction $\Phi_n$ de $\Phi$ \`a $\GG_{a,n} \croi \Omega$ donne un isomorphisme :

\[\Phi_n : \mu^* \call{F}\res{\GG_{a,n} \croi \Omega} \sta{\iso}{\to} p^*\call{F}\res{\GG_{a,n} \croi \Omega}\]
or, $p^*\call{F}\res{\GG_{a,n} \croi \Omega} \iso \kk[T]/(T^n) \tens_\kk \call{F}\res{\Omega}$.

Notons $\mu^* \sigma$ l'image de $\sigma$ dans $\Gamma(\mu\inv{\Omega},\mu^*\call{F})$ par le morphisme naturel $\call{F } \to \mu_*\mu^*\call{F}$. 

On a  pour tout $n >0$, $\Phi_n(\mu^*\sigma) \in \kk[T]/(T^n)\tens_\kk \Gamma(\Omega,\call{F})$ (et $\Phi_n(\mu^*\sigma) = \Phi_{n+1}(\mu^*\sigma) \mod (T^n)$).

Notons $d : \kk[T] \to \kk , \; P(T) \donne P'(0)$. On a $\goth g_a =\kk d$ et pour tout $n$, $\goth g_a^n = \kk d^n$ o\`u $d^n : \kk[T]  \to \kk , \; P(T) \donne P^{(n)}(0)$.

L'action de $U(\goth g_a)$ sur $\Gamma(\Omega ,\call{F})$ est telle que :
\[d^n . \sigma =(d^n \tens 1)(\Phi_n(\mu^* \sigma)) \]
pour tout $n >0$.

On a alors par hypoth\`ese :
\[\qq n \ge n_0, (d^n \tens 1)(\Phi_n(\mu^*\sigma)) = d^n . \sigma =0 \p\]

Donc il existe $\sigma_1,..., \sigma_{n_0}\in\Gamma(\Omega,\call{F})$ tels que :
\[\Phi_n(\mu^*\sigma) = 1 \tens \sigma_1 + ... + T^{n_0-1}\tens \sigma_{n_0} \mod (T^n)\]
pour tout $n \ge n_0$ (il suffit de poser $\sigma_i:=(i-1)!d^{i-1} \sigma$ pour $1 \le i \le n_0$)..

Posons : \[\sigma':= 1 \tens \sigma_1 + T \tens \sigma_2 +... + T^{n_0-1}\tens \sigma_{n_0} \in \kk[T]\tens_\kk \Gamma(\Omega, \call{F}) = \Gamma(\GG_a \croi \Omega,p^*\call{F})\p\]

Comme $\Phi_n(\mu^*\sigma) = \sigma'\res{\GG_{a,n} \croi \Omega}$ pour tout $n \ge n_0$, comme $\Phi$ est un isomorphisme de faisceaux et comme $p^*\call{F}$ est un faisceau localement libre, on a :

\[\mu^*\sigma \res{\mu\inv \Omega \cap \GG_a \croi \Omega} = \Phi\inv(\sigma')\res{\GG_a \croi \Omega \cap \mu\inv\Omega}\]

et donc (comme $\mu^* \call{F}$ est un faisceau), si l'on pose $W := \mu\inv \Omega \uni \GG_a \croi \Omega$, il existe $\chap{\sigma} \in \Gamma(W, \mu^*\call{F})$ tel que : $\mu^*\sigma = \chap{\sigma} \res{\mu\inv \Omega}$.

Soit $g \in \GG_a$. On pose $i : Y \to \GG_a \croi Y $, $y \donne (g,g\inv y)$. On a :
\[\Omega \uni g \Omega \sub i\inv W \et \mu \circ i =\Id_Y\p\]

On a donc :

\[i^*\mu^* \sigma =\sigma \; ;\]
\[i^*\chap{\sigma} \in \Gamma(i\inv W,\call{F}) \; ;\]
\[\sigma = i^*\chap{\sigma}\res{\Omega} \p\]

En particulier, on a prolong\'e $\sigma$ \`a l'ouvert $\Omega \uni g\Omega$. Puisqu'on peut le faire pour tout $g \in \GG_a$, $\sigma$ se prolonge \`a l'ouvert $\uni_{g \in \GG_a} g \Omega$ qui est stable par $\GG_a$. 

\end{dem}

\subsection{Op\'erateurs diff\'erentiels tordus}

Soit $\lambda = \lambda_1\omega_1 + \lambda_2\omega_2$ un caract\`ere de $\chap{T}$. On notera $\call{D}_\lambda$ le faisceau d'op\'erateurs diff\'erentiels sur $\adh{G}$ tordu par le faisceau inversible $\call{L}_\lambda$ :
\[\call{D}_\lambda = \call{L}_\lambda \tenso{\call{O}_{\adh{G}}} \call{D}_{\adh{G}}\tenso{\call{O}_{\adh{G}}}\call{L}_\lambda\inv \p\]

Rappelons les notations du paragraphe \ref{par:secglob} :
 \[f_\lambda := f_1^{\lambda_1} \cten f_2^{\lambda_2}\res{(\adh{G})_0} \in \Gamma ((\adh{G})_0,\call{L}_\lambda) \p\]

Notons $\Phi_\lambda$ le morphisme d'alg\`ebres :
\[\Phi_\lambda : U(\goth g \croi \goth g) \to \call{D}_\lambda(\adh{G})\]
induit par l'action de $\goth g \croi \goth g$ sur $\call{L}_\lambda$. En utilisant les formules de changement de variables de la proposition \ref{pro:cdv}, on obtient :

\begin{pro}
Dans l'anneau $\Gamma((\adh{G})_0,\call{D}_\lambda)$, on a : 
\begin{eqnarray*}
\Phi_\lambda(Y_1^{(g)}) &=& f_\lambda \tens \Phi_0(Y_1^{(g)})\tens f_\lambda\inv -\lambda_2U_{1,2} \;,\\
\Phi_\lambda(Y_2^{(g)}) &=& f_\lambda \tens \Phi_0(Y_2^{(g)})\tens f_\lambda\inv -\lambda_1 U_{2,3}\;,\\
\Phi_\lambda(Y_3^{(g)}) &=& f_\lambda \tens \Phi_0(Y_3^{(g)})\tens f_\lambda\inv-\lambda_1 U_{1,3} + \lambda_2 (U_{1,2}U_{2,3}-U_{1,3})\;,\\
\Phi_\lambda(X_1^{(g)}) &=& f_\lambda \tens \Phi_0(X_1^{(g)})\tens f_\lambda\inv\;,\\
\Phi_\lambda(X_2^{(g)}) &=& f_\lambda \tens \Phi_0(X_2^{(g)})\tens f_\lambda\inv\;,\\
\Phi_\lambda(X_3^{(g)}) &=& f_\lambda \tens \Phi_0(X_3^{(g)})\tens f_\lambda\inv\;.
\end{eqnarray*}

\end{pro}

\begin{dem}
Si $a \in \Gamma((\adh{G})_0,\call{O}_{\adh{G}})$ et si $\xi \in \goth g$, alors :
\[\xi^{(g)}.(af_\lambda) = (\xi^{(g)}.a+\frac{\xi^{(g)}.f_\lambda}{f_\lambda})f_\lambda \p\]
\end{dem}

On a aussi un op\'erateur diff\'erentiel tordu d'ordre $2$ particulier :

\begin{lem}\label{lem:globtordu}
Pour tout caract\`ere $\lambda$ de $\chap{T}$, 

\[f_\lambda \tens \partial_{\alpha_1}\partial_{\alpha_2} \tens f_\lambda\inv \in \Gamma (\adh{G},\call{D}_\lambda) \p\]
\end{lem}

\begin{dem}
Posons $D_\lambda := f_\lambda \tens \partial_{\alpha_1}\partial_{\alpha_2} \tens f_\lambda\inv $. Comme pour le th\'eor\`eme \ref{thm:opdiffglob}, il suffit de v\'erifier que $D_\lambda$ 
est $U(\goth n^- \croi n)$ et $U(\goth n \croi \goth n^-)-$fini. Puisque $D_\lambda \in \Gamma((\adh{G})_0, \call{D}_\lambda)$ et que $(\adh{G})_0$ est $B \croi B^--$stable, on sait d\'ejà que $D_\lambda$ est $U(\goth n \croi \goth n^-)$-fini. Pour le côt\'e  $U(\goth n^- \croi \goth n)-$fini, nous allons montrer que $D_\lambda \in \Gamma(B^-B, \call{D}_\lambda)$. 

Remarquons que $B^-B = \{[g] \in G \tq g_{1,1}\Delta_{3,3}\neq 0\}$.

On pose :
\[{f_\lambda}^* : (g,g') \donne (g_{1,1})^{\lambda_1} (g'_{3,3})^{\lambda_2} \]
c'est une application d\'efinie sur un ouvert de $M_3(\kk) \croi M_3(\kk)$. Si on pose $(\adh{G})_\infi := \{([g],[g']) \in \adh{G} \tq g_{1,1}g'_{3,3} \neq 0\}$, $f_\lambda^*$ d\'efinit un \'el\'ement de $\Gamma((\adh{G})_\infi,\call{L}_\lambda)$.

Dans l'espace $\Gamma(BB^- \cap B^-B,\call{L}_\lambda)$, on a l'\'egalit\'e :
\[f_\lambda = \left(\frac{g_{3,3}}{g_{1,1}}\right)^{\lambda_1} \left(\frac{g'_{1,1}}{g'_{3,3}} \right)^{\lambda_2} f_\lambda^*\]
\[
= \left(\frac{g_{3,3}}{g_{1,1}}\right)^{\lambda_1} \left(\frac{\Delta_{1,1}}{\Delta_{3,3}} \right)^{\lambda_2} f_\lambda^* \p\]

Donc en notant  $h$ la fraction rationnelle sur $G$ :
\[\left(\frac{g_{3,3}}{g_{1,1}}\right)^{-\lambda_1} \left(\frac{\Delta_{1,1}}{\Delta_{3,3}} \right)^{-\lambda_2}\]

on a :
\[D_\lambda = f_\lambda^* \tens h \inv D_0 h \tens {f_\lambda^*}\inv \p\]

Calculons dans l'anneau $\Gamma(BB^- \cap B^-B,\call{D}_{\adh{G}})$ :
\[h\inv D_0 h = h\inv \partial_{\alpha_1} \partial_{\alpha_2} h\]
\[= h\inv \partial_{\alpha_1} h \partial_{\alpha_2} + h\inv \partial_{\alpha_1}\partial_{\alpha_2}(h)\]

\[= h\inv h \partial_{\alpha_1}\partial_{\alpha_2} + h\inv \partial_{\alpha_1}(h)\partial_{\alpha_2} + h\inv\partial_{\alpha_1}(\partial_{\alpha_2}(h)) + h\inv\partial_{\alpha_2}(h)\partial_{\alpha_1}\]

\begin{eqnarray}\label{eq:e1}
= D_0 +h\inv \partial_{\alpha_1}(h)\partial_{\alpha_2} + h\inv\partial_{\alpha_2}(h)\partial_{\alpha_1} \nonumber\\
+ (h\inv \partial_{\alpha_1}(h))(h\inv\partial_{\alpha_2}(h)) + \partial_{\alpha_1}(h\inv\partial_{\alpha_2}(h))
\end{eqnarray}
Or :
\[\partial_{\alpha_1}= \frac{\Delta_{1,1}}{g_{3,3}}\partial_{g_{1,1}}\]
donc :
\[h\inv\partial_{\alpha_1}(h) = - \lambda_1\frac{\partial_{\alpha_1}(g_{3,3})}{g_{3,3}} -\lambda_2 \frac{\partial_{\alpha_1}(\Delta_{1,1})}{\Delta_{1,1}} + \frac{\partial_{\alpha_1}\left((g_{1,1})^{\lambda_1}(\Delta_{3,3})^{\lambda_2}\right)}{(g_{1,1})^{\lambda_1}(\Delta_{3,3})^{\lambda_2}} \]
\[ = \frac{\partial_{\alpha_1}((g_{1,1})^{\lambda_1}(\Delta_{3,3})^{\lambda_2})}{(g_{1,1})^{\lambda_1}(\Delta_{3,3})^{\lambda_2}} \]\begin{eqnarray}\label{eq:e2}
\in \frac{\Delta_{1,1}}{g_{3,3}}\kk[g_{i,j},(g_{1,1})\inv,(\Delta_{3,3})\inv] \p
\end{eqnarray}

De m\^eme, comme $\partial_{\alpha_2} = \frac{g_{3,3}}{\Delta_{1,1}}(\Delta_{2,2}\partial_{g_{1,1}} + \Delta_{1,1}\partial_{g_{2,2}} + \Delta_{2,1}\partial_{g_{1,2}} + \Delta_{1,2} \partial_{g_{2,1}})$, on a :
\[h\inv \partial_{\alpha_2}(h) = - \lambda_1\frac{\partial_{\alpha_2}(g_{3,3})}{g_{3,3}} -\lambda_2 \frac{\partial_{\alpha_2}(\Delta_{1,1})}{\Delta_{1,1}} + \frac{\partial_{\alpha_2}((g_{1,1})^{\lambda_1}(\Delta_{3,3})^{\lambda_2})}{(g_{1,1})^{\lambda_1}(\Delta_{3,3})^{\lambda_2}}\]
\[ = -\lambda_2\frac{g_{3,3}^2}{\Delta_{1,1}} + \frac{\partial_{\alpha_2}((g_{1,1})^{\lambda_1}(\Delta_{3,3})^{\lambda_2})}{(g_{1,1})^{\lambda_1}(\Delta_{3,3})^{\lambda_2}}\]
\begin{eqnarray}\label{eq:e3}
\in \frac{g_{3,3}}{\Delta_{1,1}}\kk[g_{i,j}, (g_{1,1})\inv,(\Delta_{3,3})\inv] \p
\end{eqnarray} 

On d\'eduit donc de (\ref{eq:e1}), (\ref{eq:e2}) et (\ref{eq:e3}) que $h\inv D_0 h $ est de la forme :
\[a_{1,1} \partial_{g_{1,1}} + a_{2,2}\partial_{g_{2,2}} + a_{1,2}\partial_{g_{1,2}} + a_{2,1}\partial_{g_{2,1}}\]
o\`u $a_{1,1},a_{2,2} ,  a_{1,2},a_{2,1} \in \kk[g_{i,j},(g_{1,1})\inv,(\Delta_{3,3})\inv]$. 

En cons\'equence : $h\inv D_0 h \in \Gamma(B^-B,\call{D}_{\adh{G}})$.
\end{dem}

\section{Irr\'eductibilit\'e des espaces de sections globales}

\begin{thm}\label{thm:irr}
Pour tout caract\`ere $\lambda$ de $\chap{T}$, le $\call{D}_\lambda(\adh{G})-$module $\Gamma(\adh{G},\call{L}_\lambda)$ est nul ou irr\'eductible.
\end{thm}

\begin{dem}
Soient $\lambda_1,\lambda_2 \in \ZZ$ tels que $\lambda = \lambda_1 \omega_1 + \lambda_2 \omega_2$ (d'o\`u : $\lambda^* = \lambda_2 \omega_1 + \lambda_1 \omega_2$).

Rappelons que \[\Gamma(\adh{G},\call{L}_\lambda) =\Plus_{m_1,m_2 \ge 0 \atop \lambda^*-m_1 \alpha_1 -m_2\alpha_2 \; \mathrm{dominant}} U(\goth g \croi \goth g). \alpha_1^{m_1}\alpha_2^{m_2} f_\lambda \p\]

On note $L(\nu)$ le $\goth g-$module simple de plus haut poids $\nu$ (pour tout caract\`ere de $\chap{T}$, $\nu$, dominant). On a un isomorphisme de $\goth g \croi \goth g-$modules irr\'eductibles :
\[U(\goth g \croi \goth g ). \alpha_1^{m_1}\alpha_2^{m_2}f_\lambda \iso \End_\kk L(\nu)\]
pour tout $\nu = \lambda^* - m_1\alpha_1 -m_2\alpha_2$ dominant. En effet, la section $\alpha_1^{m_1} \alpha_2^{m_2}f_\lambda$ est un $\cha{B}\croi \cha{B}^--$vesteur propre de poids $(\nu , -\nu)$.

Posons $\sigma_\nu :=\alpha_1^{m_1}\alpha_2^{m_2}f_\lambda \in \Gamma(BB^-,\call{L}_\lambda)$ pour tout caract\`ere $\nu=\lambda^* -m_1 \alpha_1 - m_2\alpha_2 $ de $\chap{T}$. Lorsque $m_1,m_2 \ge 0$, $\sigma_\nu \in \Gamma((\adh{G})_0,\call{L}_\lambda)$ et lorsque, de plus, $\nu$ est dominant, $\sigma_\nu \in \Gamma(\adh{G},\call{L}_\lambda)$.

Comme 
\[\Gamma(\adh{G},\call{L}_\lambda) = \Plus_\nu U(\goth g \croi \goth g).\sigma_\nu\]
(somme sur les caract\`eres $\nu$ dominants de $\lambda^* -\ZZ_{\ge 0} \alpha_1 -\ZZ_{\ge 0} \alpha_2$), comme $U(\goth g \croi \goth g)$ agit sur $\Gamma(\adh{G},\call{L}_\lambda)$ via un  morphisme d'alg\`ebres $U(\goth g \croi \goth g) \to \call{D}_\lambda(\adh{G})$ et comme les  $\goth g \croi \goth g-$modules $U(\goth g \croi \goth g).\sigma_\nu$ sont irr\'eductibles, il suffit de montrer que \[\sigma_{\nu'} \in \call{D}_\lambda(\adh{G}).\sigma_\nu\]
pour tous $\nu, \nu' \in \lambda^* - \ZZ_{\ge 0} \alpha_1 - \ZZ_{\ge 0}\alpha_2$ dominants.

Soient $m_1,m_2,m_1',m_2'$ des entiers positifs tels que $\nu = \lambda^* -m_1 \alpha_1 - m_2\alpha_2$, $\nu' = \lambda^* -m_1' \alpha_1 - m_2'\alpha_2$.

Soient $\nu_1,\nu_2 \in \ZZ_{\ge 0}$ tels que $\nu = \nu_1 \omega_1 + \nu_2 \omega_2$.

On a donc :

\begin{eqnarray}
\nu_1 = \lambda_2 - 2 m_1 + m_2 \et \nu_2 = \lambda_1 + m_1 -2 m_2 \p
\end{eqnarray}

{\bf $\mathbf 1-$er cas }: $\nu' = \nu +\alpha_1 + \alpha_2$.

Dans ce cas, $m_1 = m_1' + 1 >0$ et $m_2 = m_2' +1 >0$.

D'apr\`es le lemme \ref{lem:globtordu}, l'op\'erateur $D_\lambda := f_\lambda \tens D_0 \tens f_\lambda\inv$ est dans $\call{D}_\lambda(\adh{G})$. Or : \[D_\lambda.\sigma_\nu = (f_\lambda \tens D_0 \tens f_\lambda\inv).\alpha_1^{m_1}\alpha_2^{m_2}f_\lambda\]
\[=D_0(\alpha_1^{m_1}\alpha_2^{m_2})f_\lambda \]
\[=m_1m_2 \alpha^{m_1-1} \alpha_2^{m_2-1} f_\lambda\]
\[= m_1m_2 \sigma_{\nu'}\]
et $\sigma_{\nu'} \in \call{D}_\lambda(\adh{G}).\sigma_\nu$.
 
\vskip .5cm

Pour les cas suivants, on va utiliser des op\'erateurs diff\'erentiels de $\call{D}_\lambda(\adh{G})$ particuliers. 

Comme le faisceau $\call{L}_\lambda$ est $\cha{G}\croi \cha{G}-$lin\'earis\'e, l'alg\`ebre $\call{D}_\lambda(\adh{G})$ est un $\cha{G}\croi \cha{G}-$module. Si $w \in \cha{G}$, on note :
\[\tilda{w}:= (w,w) \in \cha{G} \croi \cha{G}\]

et on pose : $D_\lambda^w:=\tilda{w}.D_\lambda= \tilda{w}.(D_\lambda(\tilda{w}\inv.\,\cdot \,)) \in \call{D}_\lambda(\adh{G})$. 

Soit \[c:= \frac{1}{3}(H_1 +H_2) +\frac{1}{9}(H_1^2 + H_2^2+ H_1 H_2) + \frac{1}{3}(Y_1X_1 +Y_2 X_2 + Y_3 X_3) \in U(\goth g) \p\]

L'\'el\'ement $c$ est dans $Z(\goth g)$ le centre de l'alg\`ebre enveloppante $U(\goth g)$ et, pour tout caract\`ere $\mu=\mu_1\omega_1 +\mu_2 \omega_2$, agit comme une homoth\'etie sur le $\goth g-$module irr\'eductible $L(\mu)$ de plus haut poids $\mu$ :
\[\qq v \in L(\mu),\; c.v = \chi_\mu(c) v \]
o\`u $\chi_\mu (c) = \frac{\mu_1+\mu_2}{3} +\frac{\mu_1^2 +\mu_1\mu_2 +\mu_2^2}{9}$.

Identifions l'\'el\'ement $c \in Z(\goth g)$ avec l'\'el\'ement $(c,0)$ du centre $Z( \goth g \croi \goth g) $ de l'alg\`ebre enveloppante de $\goth g \croi \goth g$.

\begin{lem}\label{lem:calc}
Soit $f \in \kk[U_{i,j} \tq 1 \le i \neq j \le 3] \sub \kk[(\adh{G})_0]$. Pour tout caract\`ere $\nu$ de $\chap{T}$ qui appartient à l'ensemble :
\[\lambda^* -\ZZ_{\ge 0} \alpha_1 -\ZZ_{\ge 0} \alpha_2 \;,\]
 on a l'\'egalit\'e suivante dans $\ma \Gamma((\adh{G})_0,\call{L}_\lambda)$  :
\begin{eqnarray}
\lefteqn{(c-\chi_\nu(c))(f\sigma_\nu) }\\
&=& \frac{1}{3}  \bigg( \alpha_1 \partial_{U_{1,2}}\partial_{U_{2,1}} \\
&& + \alpha_2 (\partial_{U_{2,3}} + U_{1,2}\partial_{U_{1,3}})(\partial_{U_{3,2} + U_{2,1}\partial_{U_{3,1}}}) \\ && + \alpha_1\alpha_2 \partial_{U_{1,3}}\partial_{U_{3,1}} \bigg)(f)\sigma_\nu
\end{eqnarray}
\end{lem}
\begin{dem}
Notons $F := Z_1 \cap Z_2$ l'unique $G \croi G-$orbite ferm\'ee de $\adh{G}$ et $\call{I}_F:= \call{I}_{Z_1} + \call{I}_{Z_2}$ son id\'eal de d\'efinition dans $\call{O}_{\adh{G}}$.

On remarque que pour tout caract\`ere $\mu$ de $\cha{T}$, on a un isomorphisme de $\goth g \croi \goth g -$modules :

\begin{eqnarray}
\lefteqn{\Gamma ((\adh{G})_0,\call{L}_{\mu^*}) \; \bigm / \Gamma((\adh{G})_0, \call{L}_{\mu^*} \tens \call{I}_F)}\\
 & \iso & \Gamma((\adh{G})_0 \cap F,\call{L}_{\mu^*}\res{F})\\
&\iso & \Gamma(BB^-/B^- \croi B^-B/B , \call{L}_{G/B^- \croi G/B}(\mu,-\mu))\\
& \iso & M^*_{(\mu,-\mu)}
\end{eqnarray}
 o\`u $\call{L}_{G/B^- \croi G/B}(\mu,-\mu)$ est le faisceau inversible sur la vari\'et\'e de drapeaux $G/B^- \croi G/B$ associ\'e au caract\`ere $(\mu,-\mu)$ de $\cha{T} \croi \cha{T} $ et $M^*_{(\mu,-\mu)}$ est le dual du $\goth g \croi \goth g-$module de Verma de plus haut poids $(\mu , -\mu)$.

Tout \'el\'ement de $Z(\goth g \croi \goth g)$ agit sur $M^*_{(\mu,-\mu)}$ comme une homoth\'etie donc on a pour tout $ \sigma \in \Gamma((\adh{G})_0,\call{L}_{\nu^*})$, \[(c-\chi_\nu(c)). \sigma \in \alpha_1 \Gamma((\adh{G})_0,\call{L}_{\nu^*}) + \alpha_2\Gamma((\adh{G})_0,\call{L}_{\nu^*}) \p\]

D'un autre c\^ot\'e, gr\^ace à la proposition \ref{pro:cdv} on peut exprimer l'op\'erateur $\Phi_{\nu^*}(c) \in \Gamma(\adh{G},\call{D}_{\nu^*})$ en fonction des coordonn\'ees $\alpha_1,\alpha_2,U_{i,j},1\le i \neq j \le 3$ :
\begin{eqnarray*}
\lefteqn{\Phi_{\nu^*}(c)}\\ 
& = & \frac{1}{3} \sigma_\nu \tens \bigg( \alpha_1 \partial_{U_{1,2}}\partial_{U_{2,1}} \\
&& + \alpha_2 (\partial_{U_{2,3}} + U_{1,2}\partial_{U_{1,3}})(\partial_{U_{3,2} + U_{2,1}\partial_{U_{3,1}}}) \\ && + \alpha_1\alpha_2 \partial_{U_{1,3}}\partial_{U_{3,1}} \bigg) + ... \tens \sigma_\nu\inv 
\end{eqnarray*}

o\`u les $...$ sont mis pour des op\'erateurs dans $\Gamma((\adh{G})_0, \call{D}_{\adh{G}}) $ qui ne changent pas les degr\'es en $\alpha_1$ et en $\alpha_2$ des mon\^omes :

\[\alpha_1^{n_1}\alpha_2^{n_2}\prod_{1 \le i \neq j \le 3}U_{i,j}^{n_{i,j}} \p\]
 On a ainsi l'\'egalit\'e de l'\'enonc\'e dans :
\[ \Gamma((\adh{G})_0,\call{L}_{\nu^*}) \sub \Gamma((\adh{G})_0,\call{L}_\lambda) \p\]
\end{dem}

{\bf $\mathbf 2-$\`eme cas : $\nu' =\nu +\alpha_2$}

Dans ce cas, $m_2 =m'_2 +1 \ge 1$.

Soit $s_1 := \left(\begin{array}{ccc}
0 & 1 &0 \\
-1 & 0 &0 \\
0 & 0& 1
\end{array}
\right) \in \cha{G}$.

Gr\^ace aux formules de changement de variables de la proposition \ref{pro:cdv}, on trouve :
\begin{eqnarray*}
\tilda{s_1}\inv.\sigma_\nu &=& \frac{g_{3,3}^{\nu_2-\lambda_1}\Delta^{m_1}\Delta_{2,2}^{\nu_1}}{\Delta_{1,1}^{\lambda_2}}f_\lambda\\
& =& \alpha_1^{m_1}\alpha_2^{m_2} (\alpha_1 + U_{1,2}U_{2,1})^{\nu_1} f_\lambda \\ &=& (\alpha_1+ U_{1,2}U_{2,1})^{\nu_1} \sigma_\nu\end{eqnarray*}
\[
D_\lambda (\tilda{s_1}\inv \sigma_\nu)\] \[= m_2 \alpha_1^{m_1-1}\alpha_2^{m_2-1}(\alpha_1+U_{1,2}U_{2,1})^{\nu_1-1}( m_1(\alpha_1 + U_{1,2}U_{2,1}) +\nu_1 \alpha_1)f_\lambda\]
\[ =  m_2(\alpha_1+U_{1,2}U_{2,1})^{\nu_1-1}(m_1\sigma_{\nu + \rho} + \nu_1 \sigma_{\nu + \alpha_2})f_\lambda \]
\[
D_\lambda^{s_1}(\sigma_\nu) \] \[= \tilda{s_1} . (D_\lambda(\tilda{s_1}\inv . \sigma_\nu)) \] \[ =  m_2\alpha_1^{m_1-1}\alpha_2^{m_2-1} ( (m_1 + \nu_1) \alpha_1 + m_1 U_{1,2}U_{2,1})f_\lambda\] \[= m_2 (m_1 U_{1,2}U_{2,1}\sigma_{\nu +\rho} + (m_1+ \nu_1) \sigma_{\nu + \alpha_2})) \p 
\]

On a donc gr\^ace au lemme \ref{lem:calc} :
\begin{eqnarray}
\lefteqn{(c-\chi_{\nu + \rho}(c)). D_\lambda^{s_1}(\sigma_\nu)}\\ & = & m_2\left( \frac{m_1}{3} + (m_1+\nu_1)(\chi_{\nu +\alpha_2}(c) - \chi_{\nu +\rho}(c))\right) \sigma_{\nu +\alpha_2} \\
&=& m_2\left( \frac{m_1}{3} -(m_1+\nu_1)\frac{\nu_1 + 2}{3} \right) \sigma_{\nu +\alpha_2}\\
&=& -\frac{m_2}{3}( (m_1+\nu_1)(\nu_1+1) + \nu_1)\sigma_{\nu +\alpha_2} \p 
\end{eqnarray}

Or, comme $\nu$ est dominant, $\nu_1 \ge 0$. De plus,  $m_1\ge 0 ,m_2 \ge 1$ donc  
\[ \frac{m_2}{3}( (m_1+\nu_1)(\nu_1+1) + \nu_1) >0\]
et \[\sigma_{\nu + \alpha_2} \in \call{D}_{\lambda}(\adh{G}).\sigma_\nu \p\]

\vskip .5cm

De m\^eme on traite le cas o\`u $\nu' = \nu +\alpha_1$.

\vskip .5cm

{\bf {$ \mathbf 3-$}\`eme cas :} $\nu' = \nu -\alpha_1$.

On pose $w:=s_1s_2$ o\`u $s_2 := \left(\begin{array}{ccc}
1 & 0 & 0\\
0 & 0 & 1\\
0 & -1 & 0
\end{array}\right)$.

Comme pr\'ec\'edemment, on peut calculer $D_\lambda^w(\sigma_\nu)$. Et on trouve :
\[
(c-\chi_{\nu+ \rho}(c))(c-\chi_{\nu+ \alpha_1}(c))(c-\chi_{\nu-\alpha_1+\alpha_2 }(c))(c-\chi_{\nu+ \alpha_2}(c))(c-\chi_{\nu}(c)).D_\lambda^w \sigma_\nu \]

\[=  r\sigma_{\nu-\alpha_1}\]

o\`u :

\[r = \]
\[ -\frac{2}{3^5} (\nu_2+3)(\nu_1+m_1+1)(\nu_1+\nu_2+1)(\nu_1+\nu_2+m_2+2)(2\nu_1+\nu_2+3)\nu_1(\nu_1-1)\]

\[>0\]

et donc $\sigma_{\nu-\alpha_1} \in \call{D}_{\lambda}(\adh{G}).\sigma_\nu$.

De m\^eme, on traite le cas o\`u $\nu'=\nu -\alpha_2$.

{\bf $\mathbf 4-$\`eme cas :} $\nu' = \rho \et \nu =0$.

On a :
\begin{eqnarray}
(c-\chi_{2\rho})(c-\chi_{\rho+\alpha_1})(c-\chi_\rho)D^{w_0}_\lambda \sigma_\rho &=& -\frac{2}{3} (m_1+4)(m_2+4)\sigma_0 \p
\end{eqnarray}

Donc $\sigma_0 \in  \call{D}_\lambda({\adh{G}}).\sigma_\rho$.

{\bf Cas g\'en\'eral :} Rappelons que :
\begin{eqnarray}
\nu = \nu' + (m'_1-m_1)\alpha_1 + (m'_2-m_2)\alpha_2 \p
\end{eqnarray} 

On raisonne par r\'ecurrence sur $|m'_1-m_1| + |m'_2 - m_2|$. 

Si $m'_1=m_1$ et $m'_2=m_2$, $\sigma_{\nu'}=\sigma_\nu$ et il n'y a rien à montrer.

Si $m'_1 < m_1$, et $m'_2 < m_2$, d'apr\`es le $1-$er cas, \[\sigma_{\nu + \rho} \in \call{D}_\lambda(\adh{G}). \sigma_\nu\]
et par hypoth\`ese de r\'ecurrence,

\[\sigma_{\nu'} \in \call{D}_\lambda({\adh{G}}).\sigma_{\nu + \rho} \sub \call{D}_\lambda({\adh{G}}). \sigma_{\nu} \p\]

Si $m'_1 < m_1$ et $m'_2 \ge m_2$, $\nu + \alpha_1$ est dominant, d'apr\`es le deuxi\`eme cas, \[
\sigma_{\nu + \alpha_1} \in \call{D}_\lambda({\adh{G}}). \sigma_\nu\]
et par hypoth\`ese de r\'ecurrence :
\[\sigma_{\nu'} \in \call{D}_\lambda({\adh{G}}).\sigma_{\nu + \alpha_1} \sub \call{D}_\lambda({\adh{G}}). \sigma_{\nu} \p\]

De m\^eme si $m'_1 \ge m_1$ et $m'_2 < m_2$, $\sigma_{\nu'} \in \call{D}_\lambda({\adh{G}}). \sigma_{\nu} $.

Si $m'_1 > m_1$ et $m'_2 \ge m_2$, alors si $m'_2= m_2$, $\nu - \alpha_1$ est dominant et d'apr\`es le $3-$\`eme cas et l'hypoth\`ese de r\'ecurrence :

\[\sigma_{\nu'} \in \call{D}_\lambda({\adh{G}}).\sigma_{\nu - \alpha_1} \sub \call{D}_\lambda({\adh{G}}). \sigma_{\nu} \p\]

Si $m'_2 > m_2$, alors $ \nu - \alpha_1$ ou $\nu -\alpha_2$ est dominant et encore d'apr\`es le $3-$\`eme cas et l'hypoth\`ese de r\'ecurrence :

\[\sigma_{\nu'} \in \call{D}_\lambda({\adh{G}}).\sigma_{\nu}\]
 ou bien ni $ \nu - \alpha_1$ ni $\nu -\alpha_2$ ne sont dominants et alors n\'ecessairement : $\nu' =0$ et $\nu = \rho$ (car $\nu' < \nu$ et $\nu$ et $\nu '$ sont dominants) ; on utilise alors le $4-$i\`eme cas pour conclure.
\end{dem}

\section{Autre cas}

Notons $S_3$ l'espace des matrices $3 \croi 3$ sym\'etriques à coefficients dans $\kk$.

Soit $\adh{\cal C}$ la vari\'et\'e :

\[\adh{\cal C}:=\{([S],[S']) \in \PP(S_3) \croi \PP(S_3) \tq SS' \in \kk I_3\}\]
c'est la vari\'et\'e des \og coniques compl\`etes \fg. C'est une compactification magnifique (de rang $2$) du $\SL_3-$espace homog\`ene $\PGL_3/PSO_3$ (on peut identifier $gPSO_3 \in \PGL_3/PSO_3$ et le couple $({}^t gg, g\inv{}^tg\inv)$).

Le th\'eor\`eme \ref{thm:irr} est encore vrai si l'on remplace $\adh{G}$ par $\adh{\cal C}$.

Pour tout faisceau inversible $\call{L}$ sur $\adh{\cal C}$, soit \[\call{D}_{\call{L}} : = \call{L} \tens \call{D}_{\adh{\cal C}} \tens \call{L}\inv \p\]
\begin{thm}\label{thm:coniq}
 Pour tout faisceau inversible $\call{L}$ sur $\adh{\cal C}$, le  $\call{D}_{\call{L}}(\adh{\cal C})-$module $\Gamma(X,\call{L})$ est soit nul soit irr\'eductible.
\end{thm}

\begin{dem}
C'est la m\^eme d\'emonstration que dans le cas de $\adh{G}$. Cette fois, on utilise la grosse cellule :
\[(\adh{\cal C})_0 := \left\{([S],[S']) \in \adh{\cal C} \tq S_{1,1}S'_{3,3} \neq 0\right\} \p\]

On note pour tout $g \in \SL_3(\kk)$, pour tout $([S],[S']) \in \adh{\cal C}$, \[g.([S] ,[S']) := ([{}^tg\inv S g\inv ], [g S {}^t g])\]

et $U:=\left\{\left(\begin{array}{ccc}
1 & u_{1,2}& u_{2,3} \\
 & 1 & u_{1,3}\\
&&1
\end{array}\right) \tq \qq 1 \le i < j \le 3, u_{i,j} \in \kk \right\}$. On a un isomorphisme de vari\'et\'es alg\'ebriques :
\[U \croi \A^2 \to (\adh{\cal C})_0\]
\[(u,\left(\begin{array}{c}
x\\y
\end{array}\right)) \donne u. \left({\left[ \begin{array}{ccc}
1  &&\\
&x &\\
&&xy
\end{array}\right], \left[\begin{array}{ccc}
xy&&\\&y&\\
&&1
\end{array}\right]}\right) \p\]

On remarque alors que \[\partial_x \partial_y \in \Gamma(\adh{\cal C},\call{D}_{\adh{\cal C}}) \p\]

Cet op\'erateur joue alors le m\^eme r\^ole que $D_0$ dans le cas de $\adh{G}$.
\end{dem}

\addcontentsline{toc}{section}{R\'ef\'erences}

\bibliographystyle{plain}
\bibliography{biblio}

\end{document}